\documentclass[10pt]{amsart}
\input{epsf}
\usepackage{amssymb,latexsym}
\topskip  \numberwithin{equation}{section}
\newtheorem{theorem}{Theorem}[section]
\newtheorem{proposition}[theorem]{Proposition}
\newtheorem{corollary}[theorem]{Corollary}

\newtheorem{remark}[theorem]{Remark}

\newtheorem{example}[theorem]{Example}

\begin{document}

\pagenumbering{arabic}
\pagestyle{headings}
\def\sof{\hfill\rule{2mm}{2mm}}
\def\ls{\leq}
\def\gs{\geq}
\def\SS{\frak S}
\def\qq{{\bold q}}
\def\txx{{\frac1{2\sqrt{x}}}}
\def\tx{{\left(\txx\right)}}
\def\Bn{\mathcal{P}_n}
\def\mn{\mbox{-}}
\def\mnk{\mbox{-} k}
\def\GA{{\frak D}}
\def\GB{{\frak D}^{(2)}}
\def\AA{{\frak A}}
\def\BB{{\frak B}}
\def\CC{{\frak C}}
\def\Dp1{Dumont permutation}
\def\Dpb{Dumont permutation of the second kind}
\def\vr{\varnothing}

\title{Restricted $132$-Dumont permutations}
\maketitle

\begin{center}Toufik Mansour \footnote{Research financed by EC's
IHRP Programme, within the Research Training Network "Algebraic
Combinatorics in Europe", grant HPRN-CT-2001-00272}
\end{center}

\begin{center}{Department of Mathematics, Chalmers University of
Technology, S-41296 G\"oteborg, Sweden

        {\tt toufik@math.chalmers.se} }
\end{center}

\section*{Abstract}
A permutation $\pi$ is said to be {\em Dumont permutations of the
first kind} if each even integer in $\pi$ must be followed by a
smaller integer, and each odd integer is either followed by a
larger integer or is the last element of $\pi$ (see, for example,
\cite{Z}). In \cite{D} Dumont showed that certain classes of
permutations on $n$ letters are counted by the Genocchi numbers.
In particular, Dumont showed that the $(n+1)$st Genocchi number is
the number of Dummont permutations of the first kind on $2n$
letters.

In this paper we study the number of Dumont permutations of the
first kind on $n$ letters avoiding the pattern $132$ and avoiding (or containing exactly
once) an arbitrary pattern on $k$ letters. In several interesting
cases the generating function depends only on $k$.

\noindent{Keywords}: Dumont permutations, restricted permutations,
generating functions.
\section{Introduction}

{\bf Classical patterns}. Let $\alpha\in\SS_n$ and $\tau\in\SS_k$
be two permutations. We say that $\alpha$ {\it contains\/} $\tau$
if there exists a subsequence $1\ls i_1<i_2<\dots<i_k\ls n$ such
that $(\alpha_{i_1}, \dots,\alpha_{i_k})$ is order-isomorphic to
$\tau$; in such a context $\tau$ is usually called a {\it
pattern\/}. We say that $\alpha$ {\it avoids\/} $\tau$, or is
$\tau$-{\it avoiding\/}, if such a subsequence does not exist. The
set of all $\tau$-avoiding permutations in $\SS_n$ is denoted
$\SS_n(\tau)$. For an arbitrary finite collection of patterns $T$,
we say that $\alpha$ avoids $T$ if $\alpha$ avoids every $\tau\in
T$; the corresponding subset of $\SS_n$ is denoted $\SS_n(T)$.

While the case of permutations avoiding a single pattern has
attracted much attention, the case of multiple pattern avoidance
remains less investigated. In particular, it is natural, as the
next step, to consider permutations avoiding pairs of patterns
$\tau_1$, $\tau_2$. This problem was solved completely for
$\tau_1,\tau_2\in\SS_3$ (see \cite{SS}), and for $\tau_1\in\SS_3$
and $\tau_2\in\SS_4$ (see \cite{W}). Several recent papers
\cite{CW,Kr,MV1,MV2,MV3,MV4} deal with the case $\tau_1\in\SS_3$,
$\tau_2\in\SS_k$ for various pairs $\tau_1,\tau_2$. Another
natural question is to study permutations avoiding $\tau_1$ and
containing $\tau_2$ exactly $t$ times. Such a problem for certain
$\tau_1,\tau_2\in\SS_3$ and $t=1$ was investigated in \cite{R},
and for certain $\tau_1\in\SS_3$, $\tau_2\in\SS_k$ in
\cite{Atk1999,Kr,MV1,RWZ}. The tools involved in these papers
include generating trees, continued fractions, Chebyshev
polynomials, and Dyck words. Also, the tools involved in these
papers include many classical sequences, for example sequence of
Catalan numbers, Fibonacci numbers, and Pell numbers.

We denote the $n$th {\em Catalan number} by
$C_n=\frac{1}{n+1}\binom{2n}{n}$. The generating function for the
Catalan numbers is denoted by $C(x)$, that is, $C(x)=\sum_{n\gs0}
C_nx^n=\frac{1-\sqrt{1-4x}}{2x}$.

{\bf Generalized patterns}. In \cite{BS} generalized permutation
patterns were introduced that allow the requirement that two
adjacent letters in a pattern must be adjacent in the permutation.
We write a classical pattern with dashes between any two adjacent
letters of the pattern, say $1342$, as $1\mn3\mn4\mn2$, and if we
write, say $24\mn3\mn1$, then we mean that if this pattern occurs
in permutation $\pi\in\SS_n$, then the letters in the permutation
$\pi$ that correspond to $2$ and $4$ are adjacent (see \cite{C}).
For example, the permutation $\pi=35421$ has only two occurrences
of the pattern $23\mn1$, namely the subsequences $352$ and $351$,
whereas $\pi$ has four occurrences of the pattern $2\mn3\mn1$,
namely the subsequences $352$, $351$, $342$, and $341$.

Claesson \cite{C} presented a complete solution for the number of
permutations avoiding any single generalized pattern of length
three with exactly one adjacent pair of letters. Claesson and
Mansour \cite{CM} presented a complete solution for the number of
permutations avoiding any pair of generalized patterns of length
three with exactly one adjacent pair of letters. Kitaev \cite{Ki}
investigated simultaneous avoidance of two or more $3$-letter
generalized patterns without internal dashes. Later,
Mansour~\cite{M1,M2} (for more details see \cite{MS}) presented a
general approach to study the number of permutations avoiding
$1\mn3\mn2$ and avoiding (or containing exactly once) an arbitrary
generalized pattern.

{\bf Dumont permutations}. A permutation $\pi$ is said to be {\em
Dumont permutations of the first kind} if each even integer in
$\pi$ must be followed by a smaller integer, and each odd integer
is either followed by a larger integer or is the last element of
$\pi$ (see, for example, \cite{Z}). For example, $2143$, $3421$,
and $4213$ are all the Dumont permutations of the first kind of
length $4$.

A permutation $\pi$ is said to be {\em Dumont permutations of the
second kind} if $\pi_i<i$ for any even position $i$, and
$\pi_i\geq i$ for any odd position $i$. For example, $2143$,
$3142$, and $4132$ are all the Dumont permutations of the second
kind of length $4$.

Dumont~\cite{D} showed the number of Dumont permutations of the
first (second) kind in $\SS_{2n}$ is given by the $(n+1)$st
Genocchi number (see \cite[Sequence A001469(M3041)]{SP}).
\begin{remark}
Let $\pi\in\SS_n$ be any \Dpb; since $\pi_2<2$ we get $\pi_2=1$.
Hence, it is easy to see that there are no Dumont permutations of
the second kind in $\SS_n(132)$ for all $n\geq 4$. So, in this
paper we discuss only the case of Dumont permutations of the first
kind and refer to them simply as Dumont permutations.
\end{remark}

We define for all $r\gs2$,
\begin{equation}
Q_r(x)=1+\frac{x^2Q_{r-1}(x)}{1-x^2Q_{r-2}(x)}. \label{eqf}
\end{equation}
We denote the solution of Recurrence~\ref{eqf} with $Q_0(x)=0$ and
$Q_1(x)=1$ by $F_r(x)$, and we denote the solution of
Recurrence~\ref{eqf} with $Q_0(x)=Q_1(x)=1$ by $G_r(x)$. For
example, $F_2(x)=1+x^2$, $F_3(x)=\frac{1+x^4}{1-x^2}$,
$G_2(x)=\frac{1}{1-x^2}$, and
$G_3(x)=\frac{1-x^2+x^4}{(1-x^2)^2}$. Evidently, $F_r(x)$ and
$G_r(x)$ are rational functions in $x^2$, and for all $r\geq 1$,
\begin{equation}
F_r(x)=1+\sum_{j=1}^{r-1}
\frac{x^{2j}}{\prod_{m=r-1-j}^{r-2}(1-x^2F_m(x))}\mbox{ and }
G_r(x)=1+\sum_{j=1}^{r-1}
\frac{x^{2j}}{\prod_{m=r-1-j}^{r-2}(1-x^2G_m(x))}.
\end{equation}

\begin{example}
Using Recurrence~\ref{eqf} it is easy to see that
    $$F_4(\sqrt{x})=\sum_{n\geq0}(f_{n+2}+f_n-2)x^n
    \mbox{ and }G_4(\sqrt{x})=1+x+\sum_{n\gs2}(3\cdot2^{n-2}-1)x^n,$$
where $f_n$ is the $n$th Fibonacci number.
\end{example}

{\bf Organization of the paper}. In this paper we use generating
function techniques to study those Dumont permutations in $\SS_n$
($n\geq0$) which avoid $132$ and avoid (or contain exactly once)
an arbitrary pattern on $k$ letters. In several interesting cases
the generating function depends only on $k$.

The paper is organized as follows. The case of Dumont permutations
avoiding both $132$ and $\tau$ is treated in Section~2. We present
a simple structure for any Dumont permutation avoiding $132$. This
structure can be obtained explicitly for several interesting
cases, including classical patterns and generalized patterns. This
allows us to find explicitly some statistics on Dumont
permutations which avoid $132$. The case of avoiding $132$ and
containing another pattern $\tau$ exactly once is treated in
Section~3. Again, we find explicitly the generating function for
several interesting cases of $\tau$, including classical patterns
and generalized patterns.

Most of the explicit solutions obtained in Sections 2-4 involve
the generating functions $F_k(x)$ and $G_k(x)$.
\section{Dumont permutations which avoid $132$ and another pattern}\label{sec21}
Let $\GA$ is the set of all Dumont permutations of all sizes
including the empty permutation. Let $\GA_\tau(n)$ denote the
number of Dumont permutations in $\SS_n(132,\tau)$, and let
$\GA_\tau(x)=\sum_{n\geq 0}\GA_\tau(n)x^n$ be the corresponding
generating function. In this section we describe a method for
enumerating Dumont permutations which avoid $132$ and another
pattern and we use our method to enumerate $\GA_\tau(n)$ for
various $\tau$. We begin with an observation concerning the
structure of the Dumnot permutations of the first kind avoiding
$132$ which holds immediately from definitions.

\begin{proposition}\label{prom1}
For any $\pi\in\GA_n(132)$ such that $\pi_j=n$, there holds one of
the following assertions:
\begin{enumerate}
\item  if $n$ is odd number then $\pi=(\pi',n)$, where
$\pi'\in\GA_{n-1}(132)$;

\item  if $n$ is even number then $\pi=(\pi',n,\pi'')$ such that $\pi'$
is a \Dp1 on the numbers $n-j+1,n-j+2,\dots,n-1$, $\pi''$ is
nonempty \Dp1 on the numbers $1,2,\cdots,n-j$, and
$j=1,2,4,\ldots,n-2$ {\rm(}the minimal element of $\pi'$ cannot be
even number{\rm)}.
\end{enumerate}
\end{proposition}

\subsection{$\tau=\vr$} As a corollary of Proposition~\ref{prom1}
we find an explicit formula for the number of $132$-avoiding
Dumont permutations in $\SS_n$.

\begin{theorem}\label{th2a} The generating function for the number of
$132$-avoiding Dumont permutations in $\SS_n$ is given by
$(1+x)C(x^2)$. In other words, the number of $132$-avoiding Dumont
permutations in $\SS_n$ is given by $C_{[n/2]}$, which is the
$[n/2]$th Catalan number.
\end{theorem}
\begin{proof}
By Proposition \ref{prom1}, we have two possibilities for block
decomposition of an arbitrary $\pi\in\GA_n(132)$. Let us write an
equation for $\GA_{\vr}(x)$. The contribution of the first
decomposition above equals
$$\sum_{n\geq0}\GA_\vr(2n+1)x^{2n+1}=x\sum_{n\geq0}\GA_\vr(2n)x^{2n},$$
equivalently,
\begin{equation}
\GA_\vr(x)-\GA_\vr(-x)=x(\GA_\vr(x)+\GA_\vr(-x)).
\label{eq1}\end{equation} The contribution of the second
decomposition above equals
$$\sum_{n\geq1}\GA_\vr(2n)x^{2n}=\sum_{n\geq1}\GA_\vr(2n-1)x^{2n}
+\sum_{n\geq1}\sum_{j=0}^n\GA_\vr(2j+1)\GA_\vr(2n+2-2j)x^{2n},$$
equivalently,
\begin{equation}
\begin{array}{l}
\GA_\vr(x)+\GA_\vr(-x)-2=\\
\qquad\qquad=x(\GA_\vr(x)-\GA_\vr(-x))
+\frac{x}{2}(\GA_\vr(x)-\GA_\vr(-x))(\GA_\vr(x)+\GA_\vr(-x)-2).
\end{array}
\label{eq2}\end{equation} By putting $\GA_\vr(x)=(1+x)A(x)$ in
Equations~\ref{eq1} and \ref{eq2} it is easy to see that
$A(x)=C(x^2)$.
\end{proof}

\subsection{A classical pattern $\tau=12\ldots k$} Let us start
by the following example.

\begin{example}\label{exa1}
By definitions we have $\GA_1(x)=1$ and $\GA_{12}(x)=1+x+x^2$.
\end{example}

The case of varying $k$ is more interesting. As an extension of
Example~\ref{exa1}, let us consider the case $\tau=12\dots k$.

\begin{theorem}\label{th2b}
Let $A_k(x)=\frac{1}{2}(\GA_{12\ldots k}(x)+\GA_{12\ldots k}(-x))$
and $B_k(x)=\frac{1}{2}(\GA_{12\ldots k}(x)-\GA_{12\ldots k}(-x))$
for all $k\geq 0$. Then
$$A_k(x)=F_k(x),\qquad B_k(x)=xF_{k-1}(x),
\quad \mbox{and}\quad \GA_{12\dots k}(x)=F_k(x)+xF_{k-1}(x).$$
\end{theorem}
\begin{proof}
Using the same arguments as in the proof of Theorem~\ref{th2a} we
get
    $$\GA_{12\ldots k}(x)-\GA_{12\ldots k}(-x)=
    x(\GA_{12\ldots(k-1)}(x)+\GA_{12\ldots(k-1)}(-x)),$$
and
$$\begin{array}{l}
\GA_{12\ldots k}(x)+\GA_{12\ldots k}(-x)-2=x(\GA_{12\ldots
k}(x)-\GA_{12\ldots k} (-x))+\\
\qquad\qquad\qquad+\frac{x}{2}(\GA_{12\ldots(k-1)}(x)-\GA_{12\ldots(k-1)}(-x))(\GA_{12\ldots
k} (x)+\GA_{12\ldots k}(-x)-2).
\end{array}$$
The rest is easy to check by the definitions of $A_k$ and $B_k$.
\end{proof}

\begin{example}
Theorem~\ref{th2b}, for $k=3$, yields
$\GA_{123}(x)=\frac{1+x+x^4-x^5}{1-x^2}$. In other words, the
number of $132$-avoiding \Dp1 in $\SS_n(123)$ is given by
$1+(-1)^n$ for all $n\geq 4$, and $1$ for $n=0,1,2,3$. An another
example, Theorem~\ref{th2b}, for $k=4$, yields
$\GA_{1234}(x)=\frac{1+2x+x^2+2x^6+x^7+x^8}{(1+x)(1-x^2-x^4)}$. In
other words, the number of $132$-avoiding \Dp1 in $\SS_n(1234)$ is
$f_{n/2+2}+f_{n/2}-2$ if $n$ is even number, otherwise $2$ for all
$n\geq2$, where $f_n$ is the $n$th Fibonacci number.
\end{example}

As an extension of Theorem~\ref{th2b}, let us define
        $$\AA(x_1,x_2,x_3,\dots)=\sum_{\pi\in\GA}
        \prod_{j\gs1} x_j^{12\dots j(\pi)},$$
where $\tau(\pi)$ is the number of occurrences of $\tau$ in $\pi$.
Let
$$\begin{array}{l}
A^{(1)}(x_1,x_2,x_3,\dots)=\frac{1}{2}(\AA(x_1,x_2,x_3,\dots)+\AA(-x_1,x_2,x_3,\dots)),\\
B^{(1)}(x_1,x_2,x_3,\dots)=\frac{1}{2}(\AA(x_1,x_2,x_3,\dots)-\AA(-x_1,x_2,x_3,\dots)).
\end{array}$$
Using the same arguments as in the proof of Theorem~\ref{th2b}, we
obtain the following.

\begin{theorem}\label{th2bg} We have
$$A^{(1)}(x_1,x_2,x_3,\dots)=1+\frac{x_1^2A^{(1)}(x_1x_2,x_2x_3,x_3x_4,\dots)}{1-x_1^2x_2A^{(1)}(x_1x_2^2x_3,x_2x_3^2x_4,x_3x_4^2x_5,\dots)},$$
and
$$B^{(1)}(x_1,x_2,x_3,\dots)=x_1A^{(1)}(x_1x_2,x_2x_3,x_3x_4,\dots).$$
\end{theorem}

As an application to Theorem~\ref{th2bg}, for $x_1=x$ and $x_j=1$,
$j\gs2$, we get that
        $$B^{(1)}(x,1,1,\dots)=xA^{(1)}(x,1,1,\dots),$$
and
$$A^{(1)}(x,1,1,\dots)=\frac{1}{1-\dfrac{x^2}{1-\dfrac{x^2}{\ddots}}}=C(x^2).$$
Hence, we have $\GA_\vr(x)=(1+x)C(x^2)$ (see Theorem~\ref{th2a}).

Another application of Theorem~\ref{th2bg} is to the number of
right to left maxima. Let $\pi\in\SS_n$, $\pi_i$ is a {\em right
to left maxima} if $\pi_i>\pi_j$ for all $i<j$. We denote the
number of right to left maxima of $\pi$ by $rlm(\pi)$.
Proposition~$5$ of \cite{BCS} proved
    $$lrm(\pi)=\sum_{j\geq1} 12\dots j(\pi)(-1)^{j-1}.$$
Therefore,
$$\sum_{\pi\in\GA}x^{|\pi|}y^{rlm(\pi)}=\AA(xy,y^{-1},y,y^{-1},\dots)$$
together with Theorem~\ref{th2bg} and
$A^{(1)}(x,1,1,\dots)=C(x^2)$ we get
$$\sum_{\pi\in\GA}x^{|\pi|}y^{rlm(\pi)}=1+xC(x^2)y+\sum_{n\gs2}x^{2n-2}C^{n-1}(x^2)y^n.$$

\begin{corollary}
The generating function for the number of Dumont permutations
avoiding $132$ and having exactly $k$ right to left maxima is
given by $x^{2k-2}C^{k-1}(x^2)$ for all $k\gs2$, and $x^kC^k(x^2)$
for $k=0,1$.
\end{corollary}
\subsection{A classical pattern $\tau=2134\ldots k$} Similarly as
in Theorem~\ref{th2b}, we obtain the case $\tau=2134\ldots k$.
\begin{theorem}\label{th2c}
For all $k\gs2$,
    $$\GA_{213\dots k}(x)=G_{k-1}(x)+xG_{k-2}.$$
\end{theorem}

\begin{example}
Theorem~\ref{th2c} for $k=3,4$ yields
$\GA_{213}(x)=\frac{1+x-x^3}{1-x}$ and
$\GA_{2134}(x)=\frac{1+x-x^2-x^3+x^4}{(1-x^2)^2}$.
\end{example}
\subsection{A generalized pattern $12\mn3\mn\cdots\mn k$} In this
subsection we use the notation of generalized patterns (see
Section~1). For example, we write the classical pattern $132$ as
$1\mn3\mn2$.

By definitions, we get $\GA_{12}(x)=1+x+x^2$. So, by the same
arguments as in the proof of Theorem~\ref{th2b}, together with
            $$\GA_{12}(x)=\GA_{1\mn2}(x)=1+x+x^2,$$
we obtain the following.

\begin{theorem}\label{th2d}
For all $k\gs1$,
    $$\GA_{12\mn3\mn\ldots\mnk}(x)=\GA_{1\mn2\mn3\mn\ldots\mnk}(x)=F_k(x)+xF_{k-1}(x).$$
\end{theorem}
A comparison of Theorem~\ref{th2b} with Theorem~\ref{th2d}
suggests that there should exist a bijection between the sets
$\SS_n(1\mn3\mn2,12\mn3\mn\cdots\mnk)$ and
$\SS_n(1\mn3\mn2,1\mn2\mn3\mn\cdots\mnk)$. However, we failed to
produce such a bijection, and finding it remains a challenging
open question.

Now, let us define
        $$\BB(x_1,x_2,x_3,\dots)=\sum_{\pi\in\GA}
        x_1^{1(\pi)}\prod_{j\gs2} x_1^{12\mn3\mn\dots\mn j(\pi)},$$
where $\tau(\pi)$ is the number of occurrences of $\tau$ in $\pi$.
Let
$$\begin{array}{l}
A^{(2)}(x_1,x_2,x_3,\dots)=\frac{1}{2}(\BB(x_1,x_2,x_3,\dots)+\BB(-x_1,x_2,x_3,\dots)),\\
B^{(2)}(x_1,x_2,x_3,\dots)=\frac{1}{2}(\BB(x_1,x_2,x_3,\dots)-\BB(-x_1,x_2,x_3,\dots)).
\end{array}$$
Using the same arguments as those in the proof of
Theorem~\ref{th2b}, we get
\begin{theorem}\label{th2dg}
$$A^{(2)}(x_1,x_2,x_3,\dots)=1+\frac{x_1^2(1-x_2+x_2A^{(2)}(x_1,x_2x_3,x_3x_4,\dots))}
{1-x_1^2x_2(1-x_2x_3+x_2x_3A^{(2)}(x_1,x_2x_3^2x_4,x_3x_4^2x_5,\dots))},$$
and
$$B^{(2)}(x_1,x_2,x_3,\dots)=x_1-x_1x_2+x_1x_2A^{(2)}(x_1,x_2x_3,x_3x_4\dots).$$
\end{theorem}

Let $\pi\in\SS_n$; we say $\pi_j$ is a {\em rise} for $\pi$ if
$\pi_j<\pi_{j+1}$ for all $j=1,2,\dots,n-1$. We denote the number
of rises of $\pi$ by $rises(\pi)$. By definitions, we have
$$\sum_{\pi\in\GA}x^{|\pi|}y^{rises(\pi)}=x-xy+(1+xy)A^{(2)}(x,y,1,1,\dots),$$
so an application for Theorem~\ref{th2dg} we get

\begin{corollary}
The generating function $\sum_{\pi\in\GA}x^{|\pi|}y^{rises(\pi)}$
is given by
$$\frac{1+xy-2x^2y+2x^2y^2-(1+xy)\sqrt{1-4x^2y}}{2x^2y^2}.$$
In other words, the generating function for Dumont permutations
avoiding $1\mn3\mn2$ with exactly $k$ rises is given by $C_k
x^{2k+1}+C_{k+1}x^{2k+2}$ for all $k\geq 1$, and $1+x+x^2$ for
$k=0$, where $C_m$ is the $m$th Catalan number.
\end{corollary}

\subsection{A generalized pattern $\tau=21\mn3\mn\cdots\mn k$} In
this subsection, we use the notation of generalized patterns (see
Section~1). For example, we write the classical pattern $132$ as
$1\mn3\mn2$.

By definitions, we get $\GA_{21}(x)=1+x$. So, by the same
arguments as in the proof of Theorem~\ref{th2b} together with
            $$\GA_{21}(x)=\GA_{2\mn1}(x)=1+x,$$
we obtain the following.

\begin{theorem}\label{th2e}
For all $k\gs2$,
    $$\GA_{21\mn3\mn\ldots\mnk}(x)=\GA_{2\mn1\mn3\mn\ldots\mnk}(x)=G_{k-1}(x)+xG_{k-2}(x).$$
\end{theorem}

A comparison of Theorem~\ref{th2c} with Theorem~\ref{th2e}
suggests that there should exist a bijection between the sets
$\SS_n(1\mn3\mn2,21\mn3\mn\cdots\mnk)$ and
$\SS_n(1\mn3\mn2,2\mn1\mn3\mn\cdots\mnk)$. However, we failed to
produce such a bijection, and finding it remains a challenging
open question.

Now, let us define
        $$\CC(x_1,x_2,x_3,\dots)=\sum_{\pi\in\GA}
        x_1^{1(\pi)}\prod_{j\gs2} x_1^{21\mn3\mn\dots\mn j(\pi)},$$
where $\tau(\pi)$ is the number of occurrences of $\tau$ in $\pi$.
Let
$$\begin{array}{l}
A^{(3)}(x_1,x_2,x_3,\dots)=\frac{1}{2}(\CC(x_1,x_2,x_3,\dots)+\CC(-x_1,x_2,x_3,\dots)),\\
B^{(3)}(x_1,x_2,x_3,\dots)=\frac{1}{2}(\CC(x_1,x_2,x_3,\dots)-\CC(-x_1,x_2,x_3,\dots)).
\end{array}$$
Using the same arguments as in the proof of Theorem~\ref{th2b}, we
get the following.
\begin{theorem}\label{th2eg} We have
$$A^{(3)}(x_1,x_2,x_3,\dots)=1+\frac{x_1^2x_2A^{(3)}(x_1,x_2x_3,x_3x_4,\dots)}
{1-x_1^2x_2A^{(3)}(x_1,x_2x_3^2x_4,x_3x_4^2x_5,\dots)},$$ and
$$B^{(3)}(x_1,x_2,x_3,\dots)=x_1A^{(3)}(x_1,x_2x_3,x_3x_4\dots).$$
\end{theorem}

Let $\pi\in\SS_n$; we say that $\pi_j$ is a {\em descent} for
$\pi$ if $\pi_j>\pi_{j+1}$ for all $j=1,2,\dots,n-1$. We denote
the number of descents of $\pi$ by $descents(\pi)$. By
definitions, we have
$$\sum_{\pi\in\GA}x^{|\pi|}y^{decents(\pi)}=(1+x)A^{(3)}(x,y,1,1,\dots),$$
therefore an application for Theorem~\ref{th2eg} we get

\begin{corollary}
The generating function
$\sum_{\pi\in\GA}x^{|\pi|}y^{descents(\pi)}$ is given by
$(1+x)C(x^2y)$. In other words, the generating function for Dumont
permutations avoiding $1\mn3\mn2$ with exactly $k$ descents is
given by $C_k x^{2k+1}+C_{k}x^{2k+2}$ for all $k\geq 0$, where
$C_m$ is the $m$th Catalan number.
\end{corollary}

\subsection{A classical pattern $\tau=23\dots k1$} Again,
Proposition~\ref{prom1} gives a complete answer for $\tau=23\dots
k1$.

\begin{theorem}\label{th2f}
For all $k\geq 3$,
    $$\GA_{23\dots k1}(x)=1+x+\frac{x^2(1+x)}{1-x^2-x^2F_{k-3}(x)}.$$
\end{theorem}
\begin{proof}
Using the same arguments as in the  proof of Theorem~\ref{th2a} we
get
    $$\GA_{23\ldots k1}(x)-\GA_{23\ldots k1}(-x)=x(\GA_{23\ldots k1}(x)+\GA_{23\ldots k1}(-x)),$$
and
$$\begin{array}{l}
\GA_{23\ldots k1}(x)+\GA_{23\ldots k1}(-x)-2=x(\GA_{23\ldots
k1}(x)-\GA_{23\ldots k1} (-x))+\\
\qquad\qquad\qquad+\frac{x}{2}(\GA_{12\cdots(k-2)
}(x)-\GA_{12\ldots(k-2)}(-x))(\GA_{23\ldots k1} (x)+\GA_{23\ldots
k1}(-x)-2).
\end{array}$$
The rest is easy to check by the definitions of $F_k(x)$ together
with Theorem~\ref{th2b}.
\end{proof}

\begin{example}
Theorem~\ref{th2f}, for $k=5$, yields
$\GA_{23451}(x)=\frac{(1+x)(1-x^2-x^4)}{1-2x^2-x^4}$. In other
words, the number of \Dp1 in $\SS_n(132,23451)$ is given by
$P_{[n/2]}$, which is the $[n/2]$th Pell number for all $n\geq 2$.
\end{example}
\section{Dumont permutations which avoid $132$ and contain another pattern exactly
once} Let $\GA_{\tau;r}(n)$ denote the number of Dumont
permutations in $\SS_n(132)$ containing $\tau$ exactly $r$ times,
and let $\GA_{\tau;r}(x)=\sum_{n\geq 0}\GA_{\tau;r}(n)x^n$ be the
corresponding generating function.

\subsection{A classical pattern $\tau=12\dots k$}
\begin{theorem}\label{th3a} Let
$$A_k(x)=\frac{x^2}{1-x^2F_{k-2}(x)}A_{k-1}(x)+\frac{x^4F_{k-1}(x)}{(1-x^2F_{k-2}(x))^2}A_{k-2}(x)$$
for all $k\geq2$, where $A_1(x)=0$ and $A_2(x)=x^4$. Then for all
$k\geq2$
$$\GA_{12\dots k;1}(x)=A_k(x)+xA_{k-1}(x).$$
\end{theorem}
\begin{proof}
By Proposition~\ref{prom1}, we have two possibilities for the
block decomposition of an arbitrary $\pi$ in $\GA_n(132)$. Let us
write an equation for $\GA_{12\dots k;1}(x)$. The contribution of
the first decomposition above is
$$\sum_{n\geq0}\GA_{12\dots k;1}(2n+1)x^{2n+1}=x\sum_{n\geq0}\GA_{12\dots(k-1);1}(2n)x^{2n},$$
equivalently
\begin{equation}
\GA_{12\dots k;1}(x)-\GA_{12\dots k;1}(-x)=x(\GA_{12\dots(k-1);1}
(x)+\GA_{12\dots(k-1);1}(-x)). \label{eqaa}
\end{equation}
The contribution of the second decomposition above is
$$\begin{array}{l}
\sum\limits_{n\geq1}\GA_{12\dots k;1}(2n)x^{2n}=
\sum\limits_{n\geq1}\GA_{12\dots k;1}(2n-1)x^{2n}+\\
\qquad\qquad\qquad\qquad\qquad+\sum\limits_{n\geq1}\sum\limits_{j=0}^n\GA_{12\dots(k-1);1}(2j+1)\GA_{12\dots
k;0} (2n+2-2j)x^{2n}+\\
\qquad\qquad\qquad\qquad\qquad\qquad+\sum\limits_{n\geq1}\sum\limits_{j=0}^n\GA_{12\dots(k-1);0}(2j+1)\GA_{12\dots
k;1} (2n+2-2j)x^{2n}, \end{array}$$ equivalently
\begin{equation}
\begin{array}{l}
\GA_{12\dots k;1}(x)+\GA_{12\dots k;1}(-x)=x(\GA_{12\dots k;1} (x)-\GA_{12\dots k;1} (-x))+\\
\qquad\quad\quad\qquad+\frac{x}{2}(\GA_{12\dots(k-1);1}(x)-\GA_{12\dots(k-1);1}(-x))(\GA_{12\dots
k;0} (x)+\GA_{12\dots k;0}(-x)-2)+\\
\qquad\quad\quad\qquad\qquad+\frac{x}{2}(\GA_{12\dots(k-1);0}(x)-\GA_{12\dots(k-1);0}(-x))
(\GA_{12\dots k;1} (x)+\GA_{12\dots k;1}(-x)).\\
\end{array}
\label{eqab}\end{equation} Using Theorem~\ref{th2b},
Equation~\ref{eqaa}, Equation~\ref{eqab}, and
Definition~\ref{eqf}, we get the desired result.
\end{proof}

\begin{example}
Theorem~\ref{th3a} for $k=3$ we get
$$\GA_{123;1}(x)=\frac{x^5(1+x-x^2)}{1-x^2},$$
and for $k=4$ we get
$$\GA_{1234;1}(x)=\frac{x^7(1+x-3x^2+2x^3+3x^4+3x^5-x^6+x^7)}{(1-x^2)(1-x^2-x^4)^2}.$$
\end{example}

As an extension of Theorem~\ref{th3a}, let us consider the case
$r\geq1$. Theorem~\ref{th2bg}, for given $k$ and $r$, yields an
explicit formula for $\GA_{12\ldots k;r}(x)$. For example, for
$k=3$ and $r=0,1,2,3,4$, we have the following.

\begin{theorem}
We have

{\rm(i)} $\GA_{123;0}(x)=\dfrac{1+x+x^4-x^5}{1-x^2}$;

{\rm(ii)} $\GA_{123;1}(x)=\dfrac{x^5(1+x-x^2)}{1-x^2}$;

{\rm(iii)}
$\GA_{123;2}(x)=\dfrac{x^5(1+x^2)(1+2x-2x^2-x^3+x^4)}{(1-x^2)^2}$;

{\rm(iv)}
$\GA_{123;3}(x)=\dfrac{x^7(1+x-x^2+x^3-x^4-x^5+x^6)}{(1-x^2)^2}$;

{\rm(v)}
$\GA_{123;4}(x)=\dfrac{x^9(1+x^2)(-1-3x+3x^2+3x^3-3x^4-x^5+x^6)}{(1-x^2)^2}$.
\end{theorem}

\subsection{A classical pattern $\tau=2134\dots k$} Similarly to
Theorem~\ref{th3a}, we have

\begin{theorem}\label{th3b} Let
$$A_k(x)=\frac{x^2}{1-x^2G_{k-2}(x)}A_{k-1}(x)+\frac{x^4G_{k-1}(x)}{(1-x^2G_{k-2}(x))^2}A_{k-2}(x)$$
for all $k\geq4$, where $A_1(x)=A_2(x)=x^2$ and
$A_3(x)=\frac{x^4}{1-x^2}$. Then, for all $k\geq2$,
                $$\GA_{213\dots k;1}(x)=A_k(x)+xA_{k-1}(x).$$
\end{theorem}

\subsection{A generalized patterns $\tau=12\mn3\mn\cdots\mnk$ and
$\tau=21\mn3\mn\cdots\mnk$} Similarly to Theorem~\ref{th3a}, we
get

\begin{theorem}\label{th3c}
Let
$$A_k(x)=\frac{x^2}{1-x^2F_{k-2}(x)}A_{k-1}(x)+\frac{x^4F_{k-1}(x)}{(1-x^2F_{k-2}(x))^2}A_{k-2}(x)$$
for all $k\geq4$, where $A_1(x)=x^2$ and $A_2(x)=2x^4$. Then, for
all $k\geq2$,
                $$\GA_{12\mn3\mn\cdots\mnk;1}(x)=A_k(x)+xA_{k-1}(x).$$
\end{theorem}

As an extension of Theorem~\ref{th3c}, let us consider the case
$r\geq1$. Theorem~\ref{th2dg}, for given $k$ and $r$, yields an
explicit formula for $\GA_{12\mn3\mn\cdots\mnk;r}(x)$. For
example, for $k=3$ and $r=0,1,2,3,4$, we have the following.

\begin{theorem}
We have

{\rm (i)}  $\GA_{12\mn3;0}(x)=\dfrac{1+x+x^4-x^5}{1-x^2}$;

{\rm (ii)}
$\GA_{12\mn3;1}(x)=\dfrac{x^5(2+3x-4x^2-x^3+2x^4)}{(1-x^2)^2}$;

{\rm
(iii)}$\GA_{12\mn3;2}(x)=\dfrac{x^7(2+2x-6x^2-x^3+6x^4+x^5-2x^6)}{(1-x^2)^3}$;

{\rm (iv)}
$\GA_{12\mn3;3}(x)=\dfrac{x^7(3+5x-10x^2-9x^3+10x^4+3x^5+4x^7-5x^8-x^9+2x^10)}{(1-x^2)^4}$;

{\rm (v)}
$\GA_{12\mn3;4}(x)=\dfrac{x^9(5+5x-23x^2-7x^3+40x^4-x^5-30x^6+5x^7+5x^8-x^9+5x^{10}+x^{11}-2x^{12})}{(1-x^2)^5}$.
\end{theorem}
Similarly to Theorem~\ref{th3a}, we have

\begin{theorem}\label{th3d} Let
$$A_k(x)=\frac{x^2}{1-x^2G_{k-2}(x)}A_{k-1}(x)+\frac{x^4G_{k-1}(x)}{(1-x^2G_{k-2}(x))^2}A_{k-2}(x)$$
for all $k\geq4$, where $A_1(x)=A_2(x)=x^2$,
$A_3(x)=\frac{x^4}{1-x^2}$, and
$A_4(x)=\frac{x^6(2-x^2)}{(1-x^2)^3}$. Then, for all $k\geq2$,
                $$\GA_{21\mn3\mn\cdots\mnk;1}(x)=A_k(x)+xA_{k-1}(x).$$
\end{theorem}

As an extension of Theorem~\ref{th3d}, let us consider the case
$r\geq1$. Theorem~\ref{th2eg}, for given $k$ and $r$, yields an
explicit formula for $\GA_{21\mn3\mn\cdots\mn k;r}(x)$. For
example, for $k=3$ and $r=0,1,2,3,4$, we have the following.

\begin{theorem}
We have

{\rm (i)}  $\GA_{21\mn3;0}(x)=\dfrac{1+x+x^4-x^5}{1-x^2}$;

{\rm (ii)} $\GA_{21\mn3;1}(x)=\dfrac{x^3(1+x-x^2)}{1-x^2}$;

{\rm (iii)}
$\GA_{21\mn3;2}(x)=\dfrac{x^5(1+2x-2x^2-x^3+x^4)}{(1-x^2)^2}$;

{\rm (iv)}
$\GA_{21\mn3;3}(x)=\dfrac{x^5(1+x-x^2+x^3-x^4-x^5+x^6)}{(1-x^2)^2}$;

{\rm (v)}
$\GA_{21\mn3;4}(x)=\dfrac{x^7(1+2x-2x^2-2x^5+2x^6+x^7-x^8)}{(1-x^2)^3}$.
\end{theorem}
\section{Further results}
Here we present three different directions to generalize the
results of the previous sections. The first of these directions is
to consider one occurrence of the classical pattern $132$. For
example, the following result is true.

\begin{theorem}
There does not exist a Dumont permutation containing $132$
(classical pattern) exactly once.
\end{theorem}
\begin{proof}
Let $\pi=(\pi',n,\pi'')$ be a \Dp1 of length $n$, which contains
the pattern $132$ exactly once. It is easy to see that there does
not exist a Dumont permutation where $n=0,1,2,3$. Suppose
$n\geq4$, and let us assume by induction on $n$ that there does
not exist a \Dp1 of length $m\leq n-1$ containing $132$ exactly
once. To prove this property for $n$, let us consider the
following two cases together using Proposition~\ref{prom1}: $n$ is
either an even number, or $n$ is an odd number.
\begin{enumerate}
\item Let $n$ be an odd number. Since $\pi$ is a \Dp1, we get
$\pi''=\vr$, so $\pi$ contains $132$ exactly once if and only if
$\pi'$ contains $132$ exactly once.

\item Let $n$ be an even number. Since $\pi$ is a \Dp1 we have
$\pi''\neq\vr$. Now, let us consider two cases: either $n$ does
not appear in the occurrence of $132$, or $n$ does.
\begin{enumerate}
\item Let the occurrence of $132$ not contain the element $n$.
So, every element of $\pi'$ is greater than every element of
$\pi''$. Therefore, either $\pi'$ is a Dumont permutation of
length $m\leq n-2$ that contains $132$ exactly once, or $\pi''$ is
a \Dp1 of length $m\leq n-1$ that contains $132$ exactly once.

\item Let the occurrence of $132$ contain the element $n$. So,
$\pi=(\pi',a,n,\pi'',a+1,\pi''')$ (see~\cite{MV4}) such that
$\pi_p=n$ and $\pi_q=a+1$, where every element of $\pi'$ is
greater than every element of $\pi''$ and every element of $\pi''$
is greater than every element of $\pi'''$. Since $n$ is even
number and maximal in $\pi$ we have that $a$ is an odd number, so
$a+1$ is an even number. Therefore, by using
Proposition~\ref{prom1} we get that $p,q$ are even numbers,
$(\pi',a)$ is of odd length, and $\pi''$ is of even length. On the
other hand, $q=p+1+|\pi''|$, so $q$ is an odd number, a
contradiction.
\end{enumerate}
\end{enumerate}
Hence, by induction on $n$ we get the desired result.
\end{proof}
The second direction is to consider more than one additional
restriction. For example, the following result is true.

\begin{theorem}\label{thga}
Let $k\geq2$. The generating function for the number of Dumont
permutations in
$\SS_n(1\mn3\mn2,1\mn2\mn3\cdots\mnk,2\mn1\mn3\cdots\mnk)$ is
given by
            $$G_{k-1}(x)+xG_{k-2}(x).$$
\end{theorem}
A comparison of Theorem~\ref{thga} with Theorem~\ref{th2c}
suggests that there should exist a bijection between the sets
$\SS_n(1\mn3\mn2,2\mn1\mn3\mn\cdots\mnk)$ and
$\SS_n(1\mn3\mn2,1\mn2\mn3\cdots\mnk,2\mn1\mn3\mn\cdots\mnk)$.
However, we failed to produce such a bijection, and finding it
remains an open question.

The third direction is to consider another $3$-letter pattern
instead of $1\mn3\mn2$.

\begin{theorem} The number of \Dpb\, in $\SS_n(3\mn2\mn1)$ is the
same as the number of \Dp1 in $\SS_n(2\mn3\mn1)$ {\rm(}or in
$\SS_n(3\mn1\mn2)${\rm)} which is equal to $C_{[n/2]}$.
\end{theorem}

{\bf Acknowledgments}: The author is grateful to S. Kitaev and the anonymous refree for
their careful reading of the manuscript.

\end{document}